\documentclass[a4paper,12pt]{article}

\usepackage{amsmath}
\usepackage[xdvi]{graphicx}
\usepackage{times}
\usepackage{txfonts}

\newcommand{\bm}[1]{\mbox{\boldmath $#1$}}

\makeatletter

\@addtoreset{equation}{section}

\newcounter{def}[section]
\renewcommand{\thedef}{\stepcounter{def}\thesection.\@arabic\c@def }

\begin{document}

\begin{center}
\textbf{\LARGE{Normal Forms of $C^\infty$ Vector Fields based on the Renormalization Group}}
\end{center}

\begin{center}
Department of Applied Mathematics and Physics

Kyoto University, Kyoto, 606-8501, Japan

Hayato CHIBA \footnote{E mail address : chiba@amp.i.kyoto-u.ac.jp}
\end{center}
\begin{center}
October 20 2008
\end{center}

\begin{center}
\textbf{Abstract}
\end{center}

The normal form theory for polynomial vector fields is extended to those for $C^\infty$ 
vector fields vanishing at the origin.
Explicit formulas for the $C^\infty$ normal form and the near identity transformation
which brings a vector field into its normal form are obtained by means of the 
renormalization group method.
The dynamics of a given vector field such as the existence of invariant manifolds is investigated via its normal form.
The $C^\infty$ normal form theory is applied to prove the existence of infinitely many 
periodic orbits of two dimensional systems which is not shown from polynomial normal forms.

\section{Introduction}

The Poincar\'{e}-Durac normal form is a fundamental tool for analyzing local dynamics of vector fields
near fixed points \cite{Arn, CW, Mur}.
It gives a local coordinate change around a fixed point 
which transforms a given vector field into a simplified one in some sense. 
The normal form theory have been well developed for polynomial vector fields;
if we have a system of ordinary differential equations $dx/dt = \dot{x} = f(x)$
on $\mathbf{R}^n$ with a $C^\infty$ vector field $f$ vanishing at the origin (i.e. $f(0) = 0$),
we expand it in a formal power series as 
\begin{equation}
\dot{x} = Ax + g_2(x) + g_3(x) + \cdots ,\quad x\in \mathbf{R}^n
\end{equation}
where $A$ is a constant matrix and $g_k(x)$'s are homogeneous polynomial vector fields of degree $k$.
Then, normal forms, simplified vector fields, for polynomials $g_2, g_3, \cdots $ are calculated one after the other as summarized in Section 2.
A coordinate transformation $x \mapsto y$ which brings a given system into a normal form is of the form
\begin{equation}
x = h(y) = y + h_2(y) + h_3(y) + \cdots ,
\end{equation}
where $h_k$'s are homogeneous polynomials on $\mathbf{R}^n$ of degree $k$ that are also obtained step by step.
It is called the \textit{near identity transformation}.
Since $h(y)$ is constructed as a formal power series, it is a diffeomorphism only on a small neighborhood of the origin.
In order to investigate the local dynamics of a given system, 
usually its normal form and the near identity transformation are truncated at a finite degree.
We will refer to this method as the \textit{polynomial normal form theory}.

In this paper, we establish the \textit{$C^\infty$ normal form theory} for systems of the form
$\dot{x} = Ax + \varepsilon f(x)$ by means of the renormalization group (RG) method, where $A$ is a diagonal matrix, 
$f$ is a $C^\infty$ vector field vanishing at the origin and $\varepsilon $ is a small parameter.
The RG method has its origin in quantum field theory and was applied to
perturbation problems of differential equations by Chen, Goldenfeld and Oono \cite{CGO, CGO2}.
For a certain class of vector fields, the RG method was mathematically justified by Chiba \cite{Chi1, Chi2, Chi3}.
Our method based on the RG method allows one to calculate normal forms of vector fields without expanding in a power series.
For example if $f$ is periodic in $x$, its $C^\infty$ normal form and a near identity transformation
are also periodic. As a result, the $C^\infty$ normal form may be valid on a large open set or the 
whole phase space and it will be applicable to detect the existence of invariant manifolds of a given system.

In Sec.2, we give a brief review of the polynomial normal forms.
In Sec.3.1,  we provide a direct sum decomposition of the space of $C^\infty$ vector fields
vanishing at the origin, which extends the decomposition of polynomial vector fields used in the 
polynomial normal form theory.
Properties of the decomposition will be investigated in detail to develop the $C^\infty$ normal form theory.
In Sec.3.2, we give a definition of the $C^\infty$ normal form and explicit formulas
for calculating them are derived by means of the RG method.
In Sec.3.3, we consider the case that the linear part of a vector field is not hyperbolic.
In this case, it is proved that if a $C^\infty$ normal form has a normally hyperbolic invariant manifold $N$,
then the original system also has an invariant manifold which is diffeomorphic to $N$.
This theorem will be used to prove the existence of infinitely many periodic orbits of a two-dimensional
system in Section 4.


\section{Review of the polynomial normal forms}

In this section, we give a brief review of the polynomial normal forms
for comparison with the $C^\infty$ normal forms to be developed in the next section.
See Chow, Li and Wang \cite{CW}, Murdock \cite{Mur} for the detail.

Let us denote by $P^k(\mathbf{R}^n)$ the set of homogeneous polynomial vector fields
on $\mathbf{R}^n$ of degree $k$.
Consider the system of ordinary differential equations on $\mathbf{R}^n$
\begin{equation}
\frac{dx}{dt} = \dot{x} = Ax + \varepsilon g_2(x) + \varepsilon ^2 g_3(x) + \cdots , \quad x\in \mathbf{R}^n,
\label{2-1}
\end{equation}
where $A$ is a constant $n\times n$ matrix, $g_k\in P^k(\mathbf{R}^n)$ for $k=2,3, \cdots $,
and where $\varepsilon \in \mathbf{R}$ is a dummy parameter which is introduced to clarify 
steps of the iteration described below.
Note that if we have a system $\dot{x} = f(x)$ with the $C^\infty$ vector field $f$
satisfying $f(0) = 0$, putting $x \mapsto \varepsilon x$ and expanding the system 
$\varepsilon \dot{x} = f(\varepsilon x)$ in $\varepsilon $ yields the system (\ref{2-1}).

Let us try to simplify Eq.(\ref{2-1}) by the coordinate transformation of the form
\begin{equation}
x = y + \varepsilon h_2(y) , \quad h_2 \in P^2(\mathbf{R}^n).
\label{2-2}
\end{equation}
Substituting Eq.(\ref{2-2}) into Eq.(\ref{2-1}) provides
\begin{equation}
\left( id + \varepsilon \frac{\partial h_2}{\partial y}(y) \right) \dot{y}
 = A(y + \varepsilon h_2(y)) + \varepsilon g_2 (y + \varepsilon h_2(y)) 
       + \varepsilon ^2 g_3 (y + \varepsilon h_2(y)) + \cdots .
\label{2-3}
\end{equation}
Expanding the above in $\varepsilon $, we obtain
\begin{equation}
\dot{y} = Ay + \varepsilon \left( g_2(y) - \frac{\partial h_2}{\partial y}(y) Ay + Ah_2(y) \right)
 + \varepsilon ^2 \widetilde{g}_3(y) + \cdots ,
\label{2-4}
\end{equation}
where $\widetilde{g}_3\in P^3(\mathbf{R}^n)$.
Let us define the map $\mathcal{L}_A$ on the set of polynomial vector fields to be
\begin{equation}
\mathcal{L}_A(f)(x) = \frac{\partial f}{\partial x}(x) Ax - Af(x).
\label{2-5}
\end{equation}
Since $\mathcal{L}_A$ keeps the degree of a monomial,
it gives the linear operator from $P^k(\mathbf{R}^n)$ into $P^k(\mathbf{R}^n)$ for any integer $k$. 
Thus, the direct sum decomposition
\begin{equation}
P^k(\mathbf{R}^n) = \mathrm{Im}\, \mathcal{L}_A|_{P^k(\mathbf{R}^n)} \oplus C_k
\label{2-6}
\end{equation}
holds, where $C_k$ is a complementary subspace of $ \mathrm{Im}\, \mathcal{L}_A|_{P^k(\mathbf{R}^n)}$.
One of the convenient choices is $C_k = \mathrm{Ker}\, (\mathcal{L}_{A}|_{P^k(\mathbf{R}^n)})^*$,
where $(\mathcal{L}_{A}|_{P^k(\mathbf{R}^n)})^*$ is the adjoint operator with respect to a given inner product on $P^k(\mathbf{R}^n)$.
In particular, it is known that $(\mathcal{L}_{A}|_{P^k(\mathbf{R}^n)})^* = \mathcal{L}_{A^*}|_{P^k(\mathbf{R}^n)} $ holds for a certain inner product,
where $A^*$ denotes the adjoint matrix of $A$: 
\begin{equation}
P^k(\mathbf{R}^n) = \mathrm{Im}\, \mathcal{L}_A|_{P^k(\mathbf{R}^n)} 
                       \oplus \mathrm{Ker}\, \mathcal{L}_{A^*}|_{P^k(\mathbf{R}^n)}.
\label{2-7}
\end{equation}
Here we note that the equality $\mathcal{L}_A (f) (x)=0$ is equivalent to the equality
$f(e^{At}x) = e^{At}f(x)$ for $t\in \mathbf{R}$;
\begin{eqnarray*}
\mathrm{Ker}\, \mathcal{L}_{A^*}|_{P^k(\mathbf{R}^n)} = \{ f\in P^k(\mathbf{R}^n)\, | \, f(e^{A^*} x) = e^{A^*}f(x)\}.
\end{eqnarray*}
Since Eq.(\ref{2-4}) is written as
\begin{equation}
\dot{y} = Ay + \varepsilon (g_2 (y) - \mathcal{L}_A(h_2)(y)) + \varepsilon ^2 \widetilde{g}_3(y) + \cdots ,
\label{2-8}
\end{equation}
there exists $h_2 \in P^2(\mathbf{R}^n)$ such that 
$g_2 - \mathcal{L}_A(h_2) \in \mathrm{Ker}\, \mathcal{L}_{A^*}|_{P^2(\mathbf{R}^n)}$.

Next thing to do is to simplify $\widetilde{g}_3 \in P^3 (\mathbf{R}^n)$ by the transformation of the form
\begin{equation}
y = z + \varepsilon ^2 h_3(z),  \quad h_3\in P^3(\mathbf{R}^n).
\label{2-9}
\end{equation}
It is easy to verify that this transformation does not change the term $g_2 - \mathcal{L}_A(h_2)$
of degree two and we obtain
\begin{equation}
\dot{y} = Ay + \varepsilon (g_2 (y) - \mathcal{L}_A(h_2)(y)) 
    + \varepsilon ^2 (\widetilde{g}_3 (y) - \mathcal{L}_A(h_3)(y)) + O(\varepsilon ^3).
\label{2-10}
\end{equation}
In a similar manner to the above, we can take $h_3$ so that 
$\widetilde{g}_3 - \mathcal{L}_A(h_3) \in \mathrm{Ker}\, \mathcal{L}_{A^*}|_{P^3(\mathbf{R}^n)}$.

We proceed by induction and obtain the well-known theorem.
\\[0.2cm]
\textbf{Theorem \thedef.} \, There exists a \textit{formal} power series transformation
\begin{equation}
x = z + \varepsilon h_2(z) + \varepsilon ^2 h_3(z) + \cdots 
\label{2-11}
\end{equation}
with $h_k \in P^k(\mathbf{R}^n)$ such that Eq.(\ref{2-1}) is transformed into the system
\begin{equation}
\dot{z} = Az + \varepsilon R_2(z) + \varepsilon ^2 R_3(z) + \cdots ,
\label{2-12}
\end{equation}
satisfying $R_k \in \mathrm{Ker}\, \mathcal{L}_{A^*} \cap P^k(\mathbf{R}^n)$ for $k=2,3,\cdots $.
The transformation (\ref{2-11}) is called the \textit{near identity transformation} and
the truncated system
\begin{equation}
\dot{z} = Az + \varepsilon R_2(z) + \varepsilon ^2 R_3(z) + \cdots + \varepsilon ^m R_m(z)
\label{2-13}
\end{equation}
is called the \textit{normal form of degree $m$}.
\\[0.2cm]
\textbf{Remark \thedef.} A few remarks are in order.
The near identity transformation (\ref{2-11}) is a diffeomorphism on a small neighborhood
of the origin. Eqs.(\ref{2-11}) and (\ref{2-12}) are not convergent series in general even if Eq.(\ref{2-1})
is convergent. See Zung \cite{Zung} for the necessary and sufficient condition for the convergence
of normal forms. Note that a normal form (\ref{2-12}) is not unique.
It is because there are many different choices of $h_2$ in Eq.(\ref{2-10})
which yield the same $R_2 := g_2 - \mathcal{L}_A(h_2)$, while such different choices of $h_2$
may change $R_3, R_4, \cdots $. The simplest form among different normal forms are called the 
hyper-normal form \cite{Mur, Mur2}.

It is known that if $A = \mathrm{diag}\, (\lambda _1, \cdots  ,\lambda _n)$ is a diagonal matrix,
$\mathrm{Im}\, \mathcal{L}_A$ and $\mathrm{Ker}\, \mathcal{L}_{A^*}\, ( = \mathrm{Ker}\, \mathcal{L}_{A})$ are given by
\begin{eqnarray}
\mathrm{Im}\, \mathcal{L}_A \cap P^k(\mathbf{R}^n)
    &=& \mathrm{span} \{\, x_1^{q_1}x_2^{q_2} \cdots x_n^{q_n} \bm{e}_i
 \,\, | \,\, \sum^n_{j=1} \lambda _j q_j \neq \lambda _i ,\,\, \sum^n_{j=1}q_j = k \},\label{2-14}  \\
\mathrm{Ker}\, \mathcal{L}_{A^*} \cap P^k(\mathbf{R}^n)
  &=& \{ f \in P^k(\mathbf{R}^n) \, | \, f(e^{At}x) = e^{At}f(x)\} \nonumber \\
& = & \mathrm{span} \{ \,x_1^{q_1}x_2^{q_2} \cdots x_n^{q_n} \bm{e}_i
 \,\, | \, \sum^n_{j=1} \lambda _j q_j = \lambda _i,\,\, \sum^n_{j=1}q_j = k  \}, \label{2-15}
\end{eqnarray}
respectively, where $\bm{e}_1 , \cdots , \bm{e}_n$ are the canonical basis of $\mathbf{R}^n$.
Indeed, we can verify that
\begin{equation}
\mathcal{L}_A (x_1^{q_1}x_2^{q_2} \cdots x_n^{q_n} \bm{e}_i) = (\sum^n_{j=1} \lambda _j q_j - \lambda _i) x_1^{q_1}x_2^{q_2} \cdots x_n^{q_n}\bm{e}_i.
\label{2-16}
\end{equation}
The condition $\sum^n_{j=1} \lambda _j q_j = \lambda _i$ is called the \textit{resonance condition}.
This implies that $R_k$ consists of resonance terms of degree $k$.


\section{$C^\infty$ normal form theory}

In this section, we develop the theory of normal forms of the system
\begin{equation}
\frac{dx}{dt} = \dot{x} = Ax + \varepsilon g_2(x) + \varepsilon ^2 g_3(x) + \cdots , \quad x\in \mathbf{R}^n,
\label{3-0}
\end{equation}
for which $g_k$ is a $C^\infty$ vector field, not a polynomial in general.
We suppose that a matrix $A$ is a diagonal matrix.
If $A$ is not semi-simple, by a suitable linear transformation and the Jordan decomposition, 
we can assume that $A$ is of the form $A = \Lambda  + \varepsilon N$, where $\Lambda $ is diagonal and $N$ is nilpotent.
By replacing $g_2(x)$ to $g_2(x) + Nx$, we can assume without loss of generality that $A$ is a diagonal matrix.


\subsection{Decomposition of the space of $C^\infty$ vector fields}

Let $P_0(\mathbf{R}^n)$ be the set of polynomial vector fields on $\mathbf{R}^n$ whose degrees are 
equal to or larger than one.
Define the linear map $\mathcal{L}_A$ on $P_0(\mathbf{R}^n)$ by Eq.(\ref{2-5}).
Then, Eq.(\ref{2-7}) gives the direct sum decomposition
\begin{equation}
P_0(\mathbf{R}^n) = \mathrm{Im}\, \mathcal{L}_A \oplus \mathrm{Ker}\, \mathcal{L}_{A}.
\label{3-1}
\end{equation}
Note that $\mathrm{Ker}\, \mathcal{L}_{A^*} = \mathrm{Ker}\, \mathcal{L}_{A}$ because $A$ is diagonal by our assumption.
By the completion, the direct sum decomposition (\ref{3-1}) is extended to the set of $C^\infty$ vector fields
vanishing at the origin.
\\[0.2cm]
\textbf{Theorem \thedef.} Let $K\subset \mathbf{R}^n$ be an open set including the origin
whose closure $\bar{K}$ is compact.
Let $\mathcal{X}^\infty_0 (K)$ be the set of $C^\infty$ vector fields $f$ on $K$ 
satisfying $f(0) = 0$. Define the linear map $\mathcal{L}_A : \mathcal{X}^{\infty}_0 (K) \to \mathcal{X}^\infty_0 (K)$
by Eq.(\ref{2-5}). Then, the direct sum decomposition
\begin{equation}
\mathcal{X}^\infty_0(K) = V_I \oplus V_K
\label{3-2}
\end{equation}
holds, where
\begin{eqnarray}
& & V_I := \mathrm{Im}\, \mathcal{L}_A, \\ \label{3-3}
& & V_K := \mathrm{Ker}\, \mathcal{L}_{A}= \{ f\in \mathcal{X}^\infty_0(K) \, | \, f(e^{A}x) = e^{A}f(x)\}. \label{3-4}
\end{eqnarray}
\textbf{Proof.}
Since the set of polynomial vector fields is dense in
$\mathcal{X}^\infty_0(K)$ equipped with the $C^\infty$ topology (Hirsch \cite{Hir}),
for any $u\in \mathcal{X}^\infty_0(K)$, there exists a sequence $u_n$ in $P_0(\mathbf{R}^n)$ such that 
$u_n \to u$ as $n\to \infty$ in $\mathcal{X}^\infty_0(K)$.
Let $u_n = v_n + w_n$ with $v_n \in \mathrm{Im}\, \mathcal{L}_A|_{P_0(\mathbf{R}^n)}
 ,\, w_n \in \mathrm{Ker}\, \mathcal{L}_A|_{P_0(\mathbf{R}^n)} $
be the decomposition along the direct sum (\ref{3-1}).
Since $u_n$ is a Cauchy sequence in $\mathcal{X}^\infty_0(K)$, $u_n (x)- u_m(x)$ is sufficiently close to zero
with its derivatives uniformly on any compact subsets in $K$ if $n$ and $m$ are sufficiently large.
Hence, $u_n - u_m$ is a polynomial whose coefficients are sufficiently close to zero.
Since $v_n$ and $w_n$ consist of non-resonance and resonance terms, respectively, they do not include
common monomial vector fields.
This shows that $v_n- v_m$ and $w_n- w_m$ are also Cauchy sequences in $\mathcal{X}^\infty_0(K)$, thus 
$v_n$ and $w_n$ converge to $v$ and $w$, respectively.
Since $\mathcal{L}_A$ is a continuous operator on $\mathcal{X}^\infty_0(K)$, $\mathcal{L}_A w_n = 0$
proves $w\in \mathrm{Ker}\, \mathcal{L}_{A}$.
For $v_n\in \mathrm{Im}\, \mathcal{L}_A$, take $F_n\in P_0(\mathbf{R}^n)$ satisfying 
$v_n = \mathcal{L}_A F_n$ and $F_n \in \mathrm{Im}\, \mathcal{L}_A$ 
that is uniquely determined through Eq.(\ref{2-16});
\begin{eqnarray*}
(\mathcal{L}_A|_{\mathrm{Im}\, \mathcal{L}_A})^{-1} (x_1^{q_1}x_2^{q_2} \cdots x_n^{q_n} \bm{e}_i) 
  = (\sum^n_{j=1} \lambda _j q_j - \lambda _i)^{-1} x_1^{q_1}x_2^{q_2} \cdots x_n^{q_n}\bm{e}_i.
\end{eqnarray*}
This proves that $F_n$ is also a Cauchy sequence converging to $F\in \mathcal{X}^\infty_0(K)$ 
and $v = \mathcal{L}_A F \in\mathcal{L}_A \mathcal{X}^\infty_0(K)$.
The desired decomposition $u = v+w$ is obtained. $\Box$
\\

We define the projections $\mathcal{P}_I : \mathcal{X}^\infty_0(K) \to V_I$ 
and $\mathcal{P}_K : \mathcal{X}^\infty_0(K) \to V_K$.
For $g\in V_I$, there exists a vector field $F\in \mathcal{X}^{\infty}_0 (K)$ such that
\begin{equation}
\mathcal{L}_A(F) = \frac{\partial F}{\partial x}(x)Ax - AF(x) = g(x).
\label{3-5}
\end{equation}
Such $F(x)$ is not unique because if $F$ satisfies the above equality,
then $F + h$ with $h\in V_K$ also satisfies it.
We write $F = \mathcal{Q}(g)$ if $F$ satisfies Eq.(\ref{3-5}) and $\mathcal{P}_K(F) = 0$.
Then $\mathcal{Q}$ defines the linear map from $V_I$ to $V_I$.
In particular, we have
\begin{equation}
\mathcal{Q}\circ \mathcal{L}_A (F) = F, \quad \mathcal{L}_A \circ \mathcal{Q}(g) = g,
\label{3-7b}
\end{equation}
for any $F, g\in V_I$.
We show a few propositions which are convenient when calculating normal forms.
\\[0.2cm]
\textbf{Proposition \thedef.} The following equalities hold for any $g\in V_I$.
\begin{eqnarray}
& \mathrm{(i)} & \mathcal{P}_K \circ \mathcal{Q}(g) =0, \\
& \mathrm{(ii)} & \mathcal{Q} [Dg \cdot \mathcal{Q}(g) + D\mathcal{Q} (g) \cdot g]
          = \mathcal{P}_I[D\mathcal{Q} (g) \cdot \mathcal{Q}(g)], \\
& \mathrm{(iii)} & e^{-As} g(e^{As}x) = \frac{\partial }{\partial s}
          \left( e^{-As} \mathcal{Q}(g) (e^{As}x)\right), \quad s\in \mathbf{R},
\label{3-6}
\end{eqnarray}
where $D$ denotes the derivative with respect to $x$.
\\[0.2cm]
\textbf{Proof.} Part (i) of Prop.3.2 follows from the definition of $\mathcal{Q}$.
To prove (ii), we write $F = \mathcal{Q}(g)$.
By using Eq.(\ref{3-5}), it is easy to verify the equality
\begin{equation}
\frac{\partial }{\partial x}\left( \frac{\partial F}{\partial x}(x) F(x) \right) Ax
 - A \left( \frac{\partial F}{\partial x}(x) F(x) \right)
   = \frac{\partial g}{\partial x}(x)F(x) + \frac{\partial F}{\partial x}(x) g(x).
\label{3-7}
\end{equation}
It is rewritten as
\begin{eqnarray*}
\mathcal{L}_A[D\mathcal{Q} (g) \cdot \mathcal{Q}(g)] = Dg \cdot \mathcal{Q}(g) + D\mathcal{Q} (g) \cdot g.
\end{eqnarray*}
Applying $\mathcal{Q}$ in the both sides and using (\ref{3-7b}) proves (ii).
Part (iii) of Prop.3.2 is shown as
\begin{eqnarray*}
\frac{\partial }{\partial s} \left( e^{-As} \mathcal{Q}(g) (e^{As}x)\right) 
&= & -A e^{-As} \mathcal{Q}(g)(e^{As}x) + e^{-As} D\mathcal{Q}(g)(e^{As}x) \cdot Ae^{As}x \\
&=& e^{-As} \mathcal{L}_A \circ \mathcal{Q}(g)(e^{As}x) = e^{-As}g(e^{As}x).\quad \Box
\end{eqnarray*}

We define the Lie bracket product (commutator) $[\, \bm{\cdot}\, ,\, \bm{\cdot}\,]$ of vector fields by
\begin{equation}
[f, g](x) = \frac{\partial f}{\partial x}(x) g(x) - \frac{\partial g}{\partial x}(x) f(x).
\end{equation}
\textbf{Proposition \thedef.} If $g, h\in V_K$, 
then $\displaystyle Dg \cdot h \in V_K$ and $[g, h] \in V_K$.
\\[0.2cm]
\textbf{Proof.} \, It follows from a straightforward calculation. $\Box$
\\[0.2cm]
\textbf{Proposition \thedef.} 
For $g\in V_I$ and $h\in V_K$, the following equalities hold:
\begin{eqnarray}
& \mathrm{(i)} & \frac{\partial g}{\partial x}h \in V_I, \quad 
    \mathcal{Q}\left( \frac{\partial g}{\partial x}h\right) = \frac{\partial \mathcal{Q}(g)}{\partial x}h, \\
& \mathrm{(ii)} & \frac{\partial h}{\partial x}g \in V_I, \quad 
    \mathcal{Q}\left( \frac{\partial h}{\partial x}g\right) = \frac{\partial h}{\partial x}\mathcal{Q}(g), \\
& \mathrm{(iii)} & [g, h] \in V_I ,\quad \mathcal{Q}([g,h]) = [\mathcal{Q}(g), h]. 
\end{eqnarray}
\textbf{Proof.} \, Put $F = \mathcal{Q}(g)$. Note that $g$ and $h$ satisfy the equalities Eq.(\ref{3-5}) and
\begin{eqnarray*}
\frac{\partial h}{\partial x}(x)Ax - Ah(x) = 0.
\end{eqnarray*}
By using them, we can prove the following equalities
\begin{eqnarray}
& & \frac{\partial }{\partial x}\left( \frac{\partial F}{\partial x}(x) h(x) \right) Ax
    -A \left( \frac{\partial F}{\partial x}(x) h(x) \right) = \frac{\partial g}{\partial x}(x)h(x), \label{3-9b} \\
& & \frac{\partial }{\partial x}\left( \frac{\partial h}{\partial x}(x) F(x) \right) Ax
    -A \left( \frac{\partial h}{\partial x}(x) F(x) \right) = \frac{\partial h}{\partial x}(x)g(x),
\label{3-9}
\end{eqnarray}
which imply that $\partial g/\partial x \cdot h \in V_I$ and $\partial h/\partial x \cdot g \in V_I$.
The same calculation also shows that $\partial F/\partial x \cdot h \in V_I$ and 
 $\partial h/\partial x \cdot F \in V_I$.
Since $\mathcal{Q} = \mathcal{L}_A^{-1}$ on $V_I$, (\ref{3-9b}) and (\ref{3-9}) give (i) and (ii) of Prop.3.4, respectively.
Part (iii) immediately follows from (i) and (ii). $\Box$
\\[0.2cm]
\textbf{Remark \thedef.} Props.3.3 and 3.4 imply $[V_K, V_K]  \subset V_K $ and $[V_I, V_K] \subset V_I$.
However, $[V_I, V_I] \subset V_I$ is not true in general.


\subsection{$C^\infty$ normal forms}

Let us consider the system on $\mathbf{R}^n$ of the form
\begin{equation}
\dot{x} = Ax + \varepsilon g_1(x) + \varepsilon ^2 g_2(x) + \cdots , \quad x\in \mathbf{R}^n,
\label{3-18}
\end{equation}
where $A$ is a constant $n\times n$ diagonal matrix, $g_1 (x), \, g_2(x), \cdots \in \mathcal{X}^\infty_0(\mathbf{R}^n)$
are $C^\infty$ vector fields vanishing at the origin,
and $\varepsilon \in \mathbf{R}$ is a parameter.
To obtain a normal form of Eq.(\ref{3-18}), we use the renormalization group method.
According to \cite{Chi1}, at first, we try to construct a regular perturbation solution for Eq.(\ref{3-18}).
Put 
\begin{equation}
x = \hat{x}(t) = x_0 + \varepsilon x_1 + \varepsilon ^2 x_2 + \cdots 
\label{3-19}
\end{equation}
and substitute it into Eq.(\ref{3-18}) :
\begin{equation}
\sum^\infty_{k=0} \varepsilon ^k \dot{x}_k
 = A\sum^\infty_{k=0} \varepsilon ^k x_k + \sum^\infty_{k=1} \varepsilon ^k g_k(\sum^\infty_{j=0} \varepsilon ^j x_j). 
\label{3-20}
\end{equation}
Expanding the right hand side with respect to $\varepsilon $ and equating the coefficients of each $\varepsilon ^k$,
we obtain the system of ODEs
\begin{eqnarray}
\dot{x}_0 &=& Ax_0, \\
\dot{x}_1 &=& Ax_1 + G_1(x_0), \\
&\vdots & \nonumber \\
\dot{x}_i &=& Ax_i + G_i(x_0,x_1, \cdots  ,x_{i-1}), \label{3-21} \\
&\vdots & \nonumber
\end{eqnarray}
where the functions $G_k$ are defined through the equality
\begin{equation}
\sum^\infty_{k=1} \varepsilon ^k g_k(\sum^\infty_{j=0} \varepsilon ^j x_j)
   = \sum^\infty_{k=1}\varepsilon ^k G_k(x_0, x_1, \cdots , x_{k-1}).
\label{3-22}
\end{equation}
For example, $G_1, G_2$ and $G_3$ are given by
\begin{eqnarray}
& & G_1(x_0)= g_1(x_0), \\
& & G_2(x_0,x_1)= \frac{\partial g_1}{\partial x}(x_0)x_1 + g_2(x_0), \\
& & G_3(x_0,x_1,x_2)= \frac{1}{2}\frac{\partial ^2g_1}{\partial x^2}(x_0)x_1^2
       + \frac{\partial g_1}{\partial x}(x_0)x_2 + \frac{\partial g_2}{\partial x}(x_0)x_1 + g_3(x_0), 
\label{3-23}
\end{eqnarray}
respectively.
Since all systems are inhomogeneous linear equations, they are solved step by step.
The zeroth order equation $\dot{x}_0 = Ax_0$ is solved as $x_0(t) = e^{At}y$,
where $y\in \mathbf{R}^n$ is an initial value.
Thus, the first order equation is written as
\begin{equation}
\dot{x}_1 = Ax_1 + g_1(e^{At}y).
\label{3-28}
\end{equation}
A general solution of this system whose initial value is $x_1(0) = h^{(1)}(y)$ is given by
\begin{equation}
x_1(t) = e^{At}h^{(1)}(y) + e^{At}\int^t_{0} \! e^{-As} g_1(e^{As}y) ds. 
\label{3-29}
\end{equation}
Now we consider choosing $h^{(1)}$ so that $x_1(t)$ above takes the simplest form.
Put $\mathcal{P}_I (g_1) = g_{1I}$ and $\mathcal{P}_K (g_1) = g_{1K}$.
Then, Prop.3.2 (iii) is used to yield
\begin{eqnarray}
x_1(t) &=& e^{At}h^{(1)}(y) + e^{At}\int^t_{0} \! e^{-As} g_{1I}(e^{As}y) ds
             + e^{At}\int^t_{0} \! e^{-As} g_{1K}(e^{As}y) ds \nonumber \\
&=& e^{At}h^{(1)}(y) 
    + e^{At}\int^t_{0} \! \frac{\partial }{\partial s}\left( e^{-As} \mathcal{Q}(g_{1I})(e^{As}y) \right) ds
    + e^{At}\int^t_{0} \! g_{1K}(y) ds \nonumber \\
&=& e^{At}h^{(1)}(y) + \mathcal{Q}(g_{1I})(e^{At}y) - e^{At}\mathcal{Q}(g_{1I})(y) + e^{At}g_{1K}(y)t.
\label{3-30}
\end{eqnarray}
Putting $h^{(1)} = \mathcal{Q}(g_{1I})$, we obtain
\begin{equation}
x_1 (t) =  \mathcal{Q}(g_{1I})(e^{At}y) + g_{1K}(e^{At}y)t.
\label{3-31}
\end{equation}
Note that the term $g_{1K}(e^{At}y)t$ is so-called the \textit{secular term}.
Next thing to do is to calculate $x_2$.
A solution of the equation of $x_2$ is given by
\begin{equation}
x_2(t) = e^{At}h^{(2)}(y) + e^{At}\int^t_{0} \!e^{-As} \left( 
         \frac{\partial g_1}{\partial x}(e^{As}y) \bigl(\mathcal{Q}(g_{1I})(e^{As}y) + g_{1K}(e^{As}y)s \bigr)
           + g_2 (e^{As}y) \right) ds, 
\label{3-32}
\end{equation}
where $h^{(2)}(y) = x_2(0)$ is an initial value.
By choosing $h^{(2)}$ appropriately as above, we can show that $x_2$ is expressed as
\begin{equation}
x_2(t) =\mathcal{Q}\mathcal{P}_I (R_2) (e^{At}y)
      + \left( \mathcal{P}_K(R_2) + \frac{\partial \mathcal{Q}(g_{1I})}{\partial y} g_{1K} \right) (e^{At}y) t
 + \frac{1}{2} \frac{\partial g_{1K}}{\partial y} (e^{At}y) g_{1K} (e^{At}y) t^2,
\label{3-33}
\end{equation}
where $R_2$ is defined by
\begin{eqnarray}
R_2(y) &=& G_2(y,\mathcal{Q}(g_{1I})(y)) - \frac{\partial \mathcal{Q}(g_{1I})}{\partial y}(y)g_{1K}(y) \nonumber \\
&=& \frac{\partial g_1}{\partial y}(y) \mathcal{Q}(g_{1I})(y) + g_2(y) 
        - \frac{\partial \mathcal{Q}(g_{1I})}{\partial y}(y)g_{1K}(y).
\label{3-34}
\end{eqnarray}
These equalities are proved in Appendix with the aid of Propositions 3.2 to 3.4.
By proceeding in a similar manner, we can prove the next proposition.
\\[0.2cm]
\textbf{Proposition \thedef.} Define functions $R_k,\, k=1,2,\cdots $ on $\mathbf{R}^n$ to be 
\begin{equation}
R_1(y) = g_1(y),
\label{3-35}
\end{equation}
and
\begin{eqnarray}
R_k(y) &=& G_k(y, \mathcal{Q}\mathcal{P}_I(R_1)(y), \mathcal{Q}\mathcal{P}_I(R_2)(y),
   \cdots , \mathcal{Q}\mathcal{P}_I(R_{k-1})(y)) \nonumber \\
& & \quad \quad
  - \sum^{k-1}_{j=1} \frac{\partial \mathcal{Q}\mathcal{P}_I(R_{j})}{\partial y}(y) \mathcal{P}_K( R_{k-j})(y),
\label{3-36}
\end{eqnarray}
for $k= 2,3, \cdots $.
Then, Eq.(\ref{3-21}) has a solution
\begin{equation}
x_i = x_i(t,y) = \mathcal{Q}\mathcal{P}_I(R_i)(e^{At}y)
   + p^{(i)}_1 (e^{At}y)t + p^{(i)}_2 (e^{At}y)t^2 + \cdots  + p^{(i)}_i (e^{At}y)t^i,
\label{3-37}
\end{equation}
where $p^{(i)}_j$'s are defined by
\begin{eqnarray}
& & p^{(i)}_1(y) = \mathcal{P}_K(R_i)(y) 
     + \sum^{i-1}_{k=1}\frac{\partial \mathcal{Q}\mathcal{P}_I (R_{k})}{\partial y}(y) 
          \mathcal{P}_K(R_{i-k})(y), \label{3-38} \\
& & p^{(i)}_j(y) = \frac{1}{j}\sum^{i-1}_{k=1}\frac{\partial p^{(k)}_{j-1}}{\partial y}(y)
                     \mathcal{P}_K(R_{i-k})(y),\,\, (j= 2,3, \cdots ,i-1), \\
& & p^{(i)}_i(y) = \frac{1}{i}\frac{\partial p^{(i-1)}_{i-1}}{\partial y}(y) \mathcal{P}_K(R_1)(y), \\
& & p^{(i)}_j(y) = 0,\,\, (j>i).
\end{eqnarray}
This proposition can be proved in the same way as Prop.A.1 in Chiba \cite{Chi1},
in which Prop.3.6 is proved by induction for the case that all eigenvalues of $A$ lie on the imaginary axis.

Now we have a formal solution of Eq.(\ref{3-18}) of the form
\begin{eqnarray}
x = \hat{x}(t,y) &=& e^{At}y + \sum^\infty_{k=1} \varepsilon ^k x_k(t, y) \nonumber \\
&=& e^{At}y + \sum^\infty_{k=1} \varepsilon ^k \left( \mathcal{Q}\mathcal{P}_I(R_k) (e^{At}y)
  + p^{(k)}_1 (e^{At}y) t\right) + O(t^2).
\label{3-41}
\end{eqnarray}
This solution diverges as $t\to \infty$ because it includes polynomials in $t$.
The RG method is used to construct better approximate solutions from the above formal solution as follows 
\cite{CGO, CGO2, Chi1, Chi2, Chi3}.

We replace polynomials $t^k$ in Eq.(\ref{3-41}) by $(t-\tau)^k$,
where $\tau \in \mathbf{R}$ is a new parameter.
Next, we regard $y = y(\tau)$ as a function of $\tau$ to be determined so that we recover 
the original formal solution :
\begin{eqnarray}
\hat{x}(t,y) = e^{At}y(\tau ) + \sum^\infty_{k=1} \varepsilon ^k \left( \mathcal{Q}\mathcal{P}_I(R_k) (e^{At}y(\tau ))
  + p^{(k)}_1 (e^{At}y(\tau )) (t-\tau )\right) + O((t-\tau )^2).
\label{3-42}
\end{eqnarray}
Since $\hat{x}(t,y)$ is independent of the ``dummy" parameter $\tau$, we impose the condition
\begin{equation}
\frac{d}{d\tau}\Bigl|_{\tau = t} \hat{x}(t,y) = 0
\label{3-43}
\end{equation}
on Eq.(\ref{3-42}), which is called the RG condition. This condition provides
\begin{equation}
0 = e^{At} \frac{dy}{dt} + \sum^\infty_{k=1} \varepsilon ^k \left( 
    \frac{\partial \mathcal{Q}\mathcal{P}_I (R_k)}{\partial y} (e^{At}y) e^{At} \frac{dy}{dt}
    - p^{(k)}_1 (e^{At}y) \right) .
\label{3-44}
\end{equation}
Substituting Eq.(\ref{3-38}) yields
\begin{eqnarray}
0 &=& e^{At} \frac{dy}{dt} + \sum^\infty_{k=1} \varepsilon ^k  
    \left(\frac{\partial \mathcal{Q}\mathcal{P}_I (R_k)}{\partial y} (e^{At}y) e^{At} \frac{dy}{dt} \right) \nonumber \\
& & \quad - \sum^\infty_{k=1} \varepsilon ^k \mathcal{P}_K(R_k) (e^{At}y)
    - \sum^\infty_{k=1} \varepsilon ^k 
      \sum^{k-1}_{j=1} \frac{\partial \mathcal{Q}\mathcal{P}_I (R_j)}{\partial y} (e^{At}y)
       \mathcal{P}_K(R_{k-j}) (e^{At}y) \nonumber \\
&=&\!\! e^{At} \left( \frac{dy}{dt}\! - \!\sum^\infty_{j=1} \varepsilon ^j \mathcal{P}_K(R_j) (y) \right)
 \!+\! \sum^\infty_{k=1}\varepsilon ^k \frac{\partial \mathcal{Q}\mathcal{P}_I (R_k)}{\partial y} (e^{At}y) e^{At}
      \left( \frac{dy}{dt} \!- \!\sum^\infty_{j=1} \varepsilon ^j \mathcal{P}_K(R_j) (y) \right). \qquad
\label{3-45}
\end{eqnarray}
Now we obtain the ODE of $y$ as
\begin{equation}
\frac{dy}{dt} = \sum^\infty_{j=1} \varepsilon ^j \mathcal{P}_K(R_j) (y), 
\label{3-46}
\end{equation}
which is called the \textit{RG equation}.
Since Eq.(\ref{3-42}) is independent of $\tau$, we put $\tau = t$ to obtain
\begin{equation}
\hat{x}(t, y(t)) = e^{At} y(t) + \sum^ \infty_{j=1} \varepsilon ^j \mathcal{Q}\mathcal{P}_I(R_j)(e^{At}y(t)),
\label{3-47}
\end{equation}
where $y(t)$ is a solution of Eq.(\ref{3-46}).
This $\hat{x}(t, y(t)) $ gives an approximate solution of the system (\ref{3-18}) if the series is truncated
at some finite order of $\varepsilon $.
Since $\mathcal{P}_K(R_j)$ satisfies $\mathcal{P}_K(R_j) (e^{At}y) = e^{At} \mathcal{P}_K(R_j)(y)$,
putting $e^{At}y = z$ transforms Eqs.(\ref{3-46}) and (\ref{3-47}) into
\begin{eqnarray}
\frac{dz}{dt} &=& Az + \sum^\infty_{j=1} \varepsilon ^j \mathcal{P}_K(R_j) (z), \label{3-48}\\
\hat{x}(t, e^{-At}z(t)) &=& z(t) + \sum^ \infty_{j=1} \varepsilon ^j \mathcal{Q}\mathcal{P}_I(R_j)(z(t)),\label{3-49}
\end{eqnarray}
respectively.
Since $\mathcal{P}_K(R_j) \in V_K$, we conclude that Eqs.(\ref{3-48}) and (\ref{3-49})
give a normal form of the system (\ref{3-18}) and a near identity transformation $x\mapsto z$.
Indeed, the next theorem is reduced to Theorem 2.1 when $g_k \in P^k(\mathbf{R}^n)$.
\\[0.2cm]
\textbf{Theorem \thedef.} Define the \textit{$m$-th order near identity transformation} to be
\begin{equation}
x = z + \varepsilon \mathcal{Q}\mathcal{P}_I(R_1)(z)
   + \varepsilon ^2\mathcal{Q}\mathcal{P}_I(R_2)(z) + \cdots  + \varepsilon^m \mathcal{Q}\mathcal{P}_I(R_m)(z).
\label{3-50}
\end{equation}
Then, it transforms the system (\ref{3-18}) into the system
\begin{equation}
\dot{z} = Az + \varepsilon \mathcal{P}_K(R_1)(z)
  + \varepsilon^2 \mathcal{P}_K(R_2)(z) + \cdots + \varepsilon^m \mathcal{P}_K(R_m)(z)
  + \varepsilon ^{m+1} S(z,\varepsilon ),
\label{3-51}
\end{equation}
where $S(z,\varepsilon )$ is a $C^\infty$ function with respect to $z$ and $\varepsilon $.
We call the truncated system
\begin{equation}
\dot{z} = Az + \varepsilon \mathcal{P}_K(R_1)(z)
  + \varepsilon^2 \mathcal{P}_K(R_2)(z) + \cdots + \varepsilon^m \mathcal{P}_K(R_m)(z)
\label{3-52}
\end{equation}
the \textit{$m$-th order normal form} of Eq.(\ref{3-18}).
This system is invariant under the action of the one-parameter group $z \mapsto e^{As}z,\, s\in \mathbf{R}$.
\\[0.2cm]
\textbf{Proof.} \, By putting $z = e^{At}y$ in Eqs.(\ref{3-50}) and (\ref{3-51}), we prove that the transformation
\begin{equation}
x = e^{At}y + \varepsilon \mathcal{Q}\mathcal{P}_I(R_1)(e^{At}y)
      + \cdots  + \varepsilon^m \mathcal{Q}\mathcal{P}_I(R_m)(e^{At}y)
\label{3-53}
\end{equation}
transforms (\ref{3-18}) into the system
\begin{equation}
\dot{y} = \varepsilon \mathcal{P}_K(R_1)(y) + \cdots + \varepsilon^m \mathcal{P}_K(R_m)(y)
 + \varepsilon ^{m+1}\widetilde{S}(t, y, \varepsilon ).
\label{3-54}
\end{equation}
The proof is done by a straightforward calculation.
By substituting Eq.(\ref{3-53}) into Eq.(\ref{3-18}), the left hand side is calculated as
\begin{equation}
\frac{dx}{dt} = \left( e^{At} + \sum^m_{k=1} \varepsilon ^k 
  \frac{\partial \mathcal{Q}\mathcal{P}_I(R_k)}{\partial y} (e^{At}y) e^{At} \right) \dot{y}
  + Ae^{At}y + \sum^m_{k=1} \varepsilon ^k 
  \frac{\partial \mathcal{Q}\mathcal{P}_I(R_k)}{\partial y} (e^{At}y) Ae^{At}y.
\label{3-55}
\end{equation}
Since $\mathcal{Q}\mathcal{P}_I(R_k)$ satisfies the equality
\begin{equation}
\frac{\partial \mathcal{Q}\mathcal{P}_I(R_k)}{\partial y}(y)Ay -   A \mathcal{Q}\mathcal{P}_I(R_k)(y)
= \mathcal{L}_A \mathcal{Q}\mathcal{P}_I(R_k)(y)  = \mathcal{P}_I(R_k)(y),
\label{3-56}
\end{equation}
Eq.(\ref{3-55}) is rewritten as
\begin{eqnarray}
\frac{dx}{dt} = \left( e^{At} \!+\! \sum^m_{k=1} \varepsilon ^k 
  \frac{\partial \mathcal{Q}\mathcal{P}_I(R_k)}{\partial y} (e^{At}y) e^{At} \right) \dot{y}+ Ae^{At}y
+ \!\sum^m_{k=1} \varepsilon ^k \left( \mathcal{P}_I(R_k)(e^{At}y)
    + A\mathcal{Q}\mathcal{P}_I(R_k)(e^{At}y) \right). \quad 
\label{3-57}
\end{eqnarray}
Furthermore, $\mathcal{P}_I(R_k) = R_k - \mathcal{P}_K(R_k)$, (\ref{3-36}) and (\ref{3-57}) are put together to yield
\begin{eqnarray}
\frac{dx}{dt} &=& \left( e^{At} \!+\! \sum^m_{k=1} \varepsilon ^k 
  \frac{\partial \mathcal{Q}\mathcal{P}_I(R_k)}{\partial y} (e^{At}y) e^{At} \right) \dot{y}+ Ae^{At}y
   + \sum^m_{k=1}\varepsilon ^k A\mathcal{Q}\mathcal{P}_I(R_k)(e^{At}y) \nonumber \\
& & \quad \quad  + \sum^m_{k=1} \varepsilon ^k \Bigl( G_k(e^{At}y, \mathcal{Q}\mathcal{P}_I(R_1)(e^{At}y), 
                 \cdots , \mathcal{Q}\mathcal{P}_I(R_{k-1})(e^{At}y)) \nonumber \\
& & \quad \quad  \quad \quad - \sum^{k-1}_{j=1} \frac{\partial \mathcal{Q}\mathcal{P}_I(R_{j})}{\partial y}(e^{At}y)
               \mathcal{P}_K( R_{k-j})(e^{At}y) - \mathcal{P}_K(R_k) (e^{At}y) \Bigr) .
\label{3-58}
\end{eqnarray}
On the other hand, the right hand side of Eq.(\ref{3-18}) is transformed as
\begin{eqnarray}
& & A (e^{At}y  + \sum^m_{k=1} \varepsilon ^k \mathcal{Q}\mathcal{P}_I(R_k) (e^{At}y))
 +  \sum^\infty_{k=1} \varepsilon ^k g_k(e^{At}y 
     + \sum^m_{j=1} \varepsilon ^j \mathcal{Q}\mathcal{P}_I(R_j) (e^{At}y)) \nonumber \\
&=& Ae^{At}y  + \sum^m_{k=1} \varepsilon ^k A\mathcal{Q}\mathcal{P}_I(R_k) (e^{At}y) \nonumber \\
& & \quad  + \sum^m_{k=1} \varepsilon ^k G_k(e^{At}y, \mathcal{Q}\mathcal{P}_I(R_1)(e^{At}y), 
                 \cdots , \mathcal{Q}\mathcal{P}_I(R_{k-1})(e^{At}y)) + O(\varepsilon ^{m+1}).
\label{3-59}
\end{eqnarray}
Thus Eq.(\ref{3-18}) is transformed into the system
\begin{eqnarray*}
\dot{y} &=& \left( e^{At} \!+\! \sum^m_{k=1} \varepsilon ^k 
  \frac{\partial \mathcal{Q}\mathcal{P}_I(R_k)}{\partial y} (e^{At}y) e^{At} \right)^{-1} \times \nonumber \\
 & & \quad        \sum^m_{k=1} \varepsilon ^k \left( \mathcal{P}_K(R_k) (e^{At}y) 
     + \sum^{k-1}_{j=1} \frac{\partial \mathcal{Q}\mathcal{P}_I(R_{j})}{\partial y}(e^{At}y)
     \mathcal{P}_K( R_{k-j})(e^{At}y) \right) + O(\varepsilon ^{m+1}) \nonumber \\
&=& e^{-At} \left(id + \sum^\infty_{j=1}(-1)^j \left( \sum^m_{k=1} \varepsilon ^k 
  \frac{\partial \mathcal{Q}\mathcal{P}_I(R_k)}{\partial y} (e^{At}y) \right)^j \right) \times \nonumber \\
& & \quad \left( e^{At} \sum^m_{i=1} \varepsilon ^i \mathcal{P}_K(R_i) (y) 
  + \sum^m_{k=1} \varepsilon ^k \frac{\partial \mathcal{Q}\mathcal{P}_I(R_{k})}{\partial y}(e^{At}y)e^{At}
     \sum^{m-k}_{i=1}\varepsilon ^i\mathcal{P}_K( R_{i})(y)  \right) + O(\varepsilon ^{m+1}) \nonumber \\
&=& \sum^m_{k=1} \varepsilon ^k \mathcal{P}_K(R_k) (y) 
  + e^{-At} \sum^\infty_{j=1} (-1)^j \left( \sum^m_{k=1} \varepsilon ^k 
  \frac{\partial \mathcal{Q}\mathcal{P}_I(R_k)}{\partial y} (e^{At}y) \right)^j
  e^{At} \sum^m_{i=m-k+1} \varepsilon ^i \mathcal{P}_K(R_i) (y) + O(\varepsilon ^{m+1}) \nonumber \\
&=& \sum^m_{k=1} \varepsilon ^k \mathcal{P}_K(R_k) (y) + O(\varepsilon ^{m+1}).
\end{eqnarray*}
This proves that Eq.(\ref{3-18}) is transformed into the system Eq.(\ref{3-54}). $\Box$
\\[0.2cm]
\textbf{Remark \thedef.} Eq.(\ref{3-51}) is valid on a region including the origin on which the near identity transformation
(\ref{3-50}) is a diffeomorphism.
In the polynomial normal form theory described in Section 2,
since $\varepsilon ^k\mathcal{Q}\mathcal{P}_I(R_k)(z)$ is a polynomial in $z$ of degree $k$,
the near identity transformation may not be a diffeomorphism when $z\sim O(1/\varepsilon )$ in general.
For the $C^\infty$ normal form, the near identity transformation may be a diffeomorphism
on larger set. For example if $\mathcal{Q}\mathcal{P}_I(R_k)(z),\, k=1,2,\cdots ,m$ are periodic 
as Example 4.1 below, Eq.(\ref{3-50}) is a diffeomorphism for any $z\in \mathbf{R}^n$ if $\varepsilon $
is sufficiently small.


\subsection{Non-hyperbolic case}

If the matrix $A$ in Eq.(\ref{3-18}) is hyperbolic, which means that no eigenvalues of $A$
lie on the imaginary axis, then the flow of Eq.(\ref{3-18}) near the origin is topologically 
conjugate to the linear system $\dot{x} = Ax$ and the local stability of the origin is easily understood.
If $A$ has eigenvalues on the imaginary axis, Eq.(\ref{3-18}) has a center manifold at the origin
and nontrivial phenomena, such as bifurcations, may occur on the center manifold.
We consider such a situation in this subsection.
By using the center manifold reduction \cite{Carr, Chi2}, we assume that all eigenvalues of $A$ lie on the 
imaginary axis. We also suppose that $A$ is diagonal as before.
In this case, the operators $\mathcal{P}_K$ and $\mathcal{Q}\mathcal{P}_I$ are calculated as follows:

Recall that the equality 
\begin{eqnarray}
\int^t_{0} \! e^{-A(s-t)}g(e^{A(s-t)}x)ds
 &=& \int^t_{0} \!  e^{-A(s-t)}\mathcal{P}_I(g)(e^{A(s-t)}x) ds
       + \int^t_{0} \!  e^{-A(s-t)}\mathcal{P}_K(g)(e^{A(s-t)}x) ds \nonumber \\
&=& \mathcal{Q}\mathcal{P}_I (g)(x) - e^{At}\mathcal{Q}\mathcal{P}_I (g)(e^{-At}x) + \mathcal{P}_K(g)(x)t
\label{hyp1}
\end{eqnarray}
holds.
We have to calculate $\mathcal{Q}\mathcal{P}_I (g)$ and $\mathcal{P}_K(g)$ to obtain the normal form (\ref{3-52}).
Since $e^{-As}g(e^{As}x)$ is an almost periodic function with respect to $s$,
it is expanded in a Fourier series as 
$e^{-As}g(e^{As}x) = \sum_{\lambda _i \in \Lambda } c(\lambda _i, x)e^{\sqrt{-1} \lambda _i s}$,
where $\Lambda$ is the set of the Fourier exponents and $c(\lambda _i, x)\in \mathbf{R}^n$
is a Fourier coefficient.
In particular, the Fourier coefficient $c(0, x)$ associated with the zero Fourier exponent
is the average of $e^{-As}g(e^{As}x)$:
\begin{equation}
c(0, x) = \lim_{t\to \infty} \frac{1}{t} \int^t \! e^{-As}g(e^{As}x) ds. 
\end{equation} 
Thus we obtain 
\begin{eqnarray}
\int^t_{0} \! e^{-A(s-t)}g(e^{A(s-t)}x)ds
 &=& \int^t_{0} \!  \sum_{\lambda _i \in \Lambda } c(\lambda _i, x)e^{\sqrt{-1} \lambda _i (s-t)} ds \nonumber \\
&=& \sum_{\lambda _i \neq 0}\frac{1}{\sqrt{-1}\lambda _i} c(\lambda _i, x) (1- e^{-\sqrt{-1}\lambda _i t})
 + c(0,x)t.
\label{hyp2}
\end{eqnarray}
Comparing it with Eq.(\ref{hyp1}), we obtain
\begin{eqnarray}
& & \mathcal{P}_K(g)(x) = c(0,x) = \lim_{t\to \infty} \frac{1}{t} \int^t \! e^{-As}g(e^{As}x) ds, \label{hyp3-1}\\
& & \mathcal{Q}\mathcal{P}_I(g)(x) 
   = \sum_{\lambda _i \neq 0}\frac{1}{\sqrt{-1}\lambda _i} c(\lambda _i, x)
   = \lim_{t\to 0} \int^t \! \left( e^{-As}g(e^{As}x) - \mathcal{P}_K(g)(x) \right)ds, \label{hyp3-2}
\end{eqnarray}
where $\int^t $ denotes the indefinite integral whose integral constant is chosen to be zero.
These formulas for $\mathcal{P}_K$ and $\mathcal{Q}\mathcal{P}_I$ allow one to calculate
the normal forms systematically.

Now we suppose that the normal form for Eq.(\ref{3-18}) satisfies
$\mathcal{P}_K(R_1) = \cdots  = \mathcal{P}_K(R_{m-1}) = 0$ for some integer $m \geq 1$.
By putting $z = e^{At}y$, Eq.(\ref{3-51}) takes the form
\begin{equation}
\dot{y} = \varepsilon ^m \mathcal{P}_K(R_{m}) + O(\varepsilon ^{m+1}).
\label{hyp4}
\end{equation}
If $\varepsilon $ is sufficiently small, some properties of Eq.(\ref{hyp4}) are obtained from
the truncated system $\dot{y} = \varepsilon ^m \mathcal{P}_K(R_{m})$.
In this manner, we can prove the next theorem.
\\[0.2cm]
\textbf{Theorem \thedef \, \cite{Chi1, Chi3}.}  Suppose that all eigenvalues of the diagonal matrix $A$ lie on the 
imaginary axis and that the normal form for Eq.(\ref{3-18}) satisfies
$\mathcal{P}_K(R_1) = \cdots  = \mathcal{P}_K(R_{m-1}) = 0$ and $\mathcal{P}_K (R_m) \neq 0$ for some integer $m \geq 1$.
If the truncated system $dy/dt = \varepsilon ^m \mathcal{P}_K(R_m)(y)$ has a normally hyperbolic invariant
manifold $N$, then for sufficiently small $|\varepsilon |$, the system (\ref{3-18}) has an 
invariant manifold $N_\varepsilon $, which is diffeomorphic to $N$.
In particular the stability of $N_\varepsilon $ coincides with that of $N$.
\\[-0.2cm]

This theorem is proved in Chiba \cite{Chi3} in terms of the RG method and a perturbation theory of invariant manifolds \cite{Wig}.
For many examples, $m=1$ and thus the dynamics of the original system (\ref{3-18}) is investigated via the first order normal form
\begin{equation}
\frac{dy}{dt} = \varepsilon \mathcal{P}_K(R_1)(y) = \varepsilon \mathcal{P}_K(g_1)(y)
 = \varepsilon \cdot \lim_{t\to \infty} \frac{1}{t}\int^t e^{-As}g_1(e^{As}y)ds,
\end{equation}
which recovers the classical averaging method.
See \cite{Chi3, Chi4} for many applications for the degenerate cases $m\geq 2$ and relationships with other perturbation methods.


\section{Examples}

In this section, we give a few examples to demonstrate our theorems.
\\[0.2cm]
\textbf{Example \thedef.} Consider the system on $\mathbf{R}^2$
\begin{equation}
\left\{ \begin{array}{l}
\dot{x}_1 = x_2 + 2\varepsilon \sin x_1, \\
\dot{x}_2 = -x_1, \\
\end{array} \right.
\label{4-1}
\end{equation}
where $\varepsilon >0$ is a small parameter.
We put $x_1 = z_1 + z_2,\, x_2 = i(z_1 - z_2)$ to diagonalize Eq.(\ref{4-1}) as
\begin{equation}
\frac{d}{dt}\left(
\begin{array}{@{\,}c@{\,}}
z_1 \\
z_2
\end{array}
\right) = \left(
\begin{array}{@{\,}cc@{\,}}
i & 0 \\
0 & -i
\end{array}
\right) \left(
\begin{array}{@{\,}c@{\,}}
z_1 \\
z_2
\end{array}
\right) + \varepsilon \left(
\begin{array}{@{\,}c@{\,}}
\sin (z_1 + z_2) \\
\sin (z_1 + z_2)
\end{array}
\right) ,
\label{4-2}
\end{equation}
where $i = \sqrt{-1}$.
We calculate the normal forms of this system in two different ways, the polynomial normal form
and the $C^\infty$ normal form.
\\[0.2cm]
\textbf{(I)} \, To calculate the polynomial normal form, we expand $\sin (z_1 + z_2)$ as
\begin{equation}
\frac{d}{dt}\left(
\begin{array}{@{\,}c@{\,}}
z_1 \\
z_2
\end{array}
\right) =  \left(
\begin{array}{@{\,}c@{\,}}
iz_1 \\
-iz_2
\end{array}
\right) + \varepsilon \left(
\begin{array}{@{\,}c@{\,}}
z_1 + z_2 \\
z_1 + z_2
\end{array}
\right) - \frac{\varepsilon }{6} \left(
\begin{array}{@{\,}c@{\,}}
(z_1 + z_2)^3 \\
(z_1 + z_2)^3
\end{array}
\right) + \frac{\varepsilon }{120} \left(
\begin{array}{@{\,}c@{\,}}
(z_1 + z_2)^5 \\
(z_1 + z_2)^5
\end{array}
\right) - \frac{\varepsilon }{5040} \left(
\begin{array}{@{\,}c@{\,}}
(z_1 + z_2)^7 \\
(z_1 + z_2)^7
\end{array}
\right) + \cdots .
\label{4-3}
\end{equation}
The fourth order normal form of this system is given by
\begin{eqnarray}
\frac{d}{dt} \left(
\begin{array}{@{\,}c@{\,}}
y_1 \\
y_2
\end{array}
\right) &=& \left(
\begin{array}{@{\,}c@{\,}}
i y_1 \\
-i y_2
\end{array}
\right) + \varepsilon \left(
\begin{array}{@{\,}c@{\,}}
y_1 \\
y_2
\end{array}
\right) - \frac{\varepsilon }{2}\left(
\begin{array}{@{\,}c@{\,}}
y_1(y_1y_2 + i) \\
y_2(y_1y_2 - i)
\end{array}
\right) \nonumber \\
& & +\frac{\varepsilon }{12}\left(
\begin{array}{@{\,}c@{\,}}
y_1^2y_2 (y_1 y_2 + 6i) \\
y_1y_2^2 (y_1 y_2 - 6i)
\end{array}
\right) - \frac{\varepsilon }{144}
\left(
\begin{array}{@{\,}c@{\,}}
y_1 (y_1^3 y_2^3 + 39i y_1^2y_2^2 + 54y_1y_2 + 18i) \\
y_2 (y_1^3 y_2^3 - 39i y_1^2y_2^2 + 54y_1y_2 - 18i) 
\end{array}
\right) .
\label{4-4}
\end{eqnarray}
Putting $y_1 = re^{i\theta },\, y_2 = re^{-i\theta } $ yields
\begin{equation}
\left\{ \begin{array}{l}
\displaystyle \dot{r} = \varepsilon r - \frac{\varepsilon}{2}r^3 + \frac{\varepsilon }{12}r^5
 - \frac{\varepsilon }{144} (r^7 + 54 r^3), \\
\displaystyle \dot{\theta }=
 1 - \frac{\varepsilon }{2} + \frac{\varepsilon }{12}\cdot 6r^2 
  - \frac{\varepsilon }{144} (39 r^4 + 18) .\\
\end{array} \right.
\label{4-5}
\end{equation}
Fixed points of the equation of $r$ (i.e. the zeros of the right hand side) 
imply periodic orbits of the original system (\ref{4-1}).
The near identity transformation is given by
\begin{equation}
\left(
\begin{array}{@{\,}c@{\,}}
z_1\\
z_2
\end{array}
\right) = \left(
\begin{array}{@{\,}c@{\,}}
y_1 \\
y_2
\end{array}
\right) + \varepsilon i \left(
\begin{array}{@{\,}c@{\,}}
\displaystyle \frac{1}{2}y_2 + \frac{1}{24}(2y_1^3 - 6y_1y_2^2 - y_2^3) + O(y_1^5, y_2^5) \\
\displaystyle -\frac{1}{2}y_1 + \frac{1}{24}(y_1^3 + 6y_1^2y_2 - 2y_2^3) + O(y_1^5, y_2^5)
\end{array}
\right),
\label{4-6}
\end{equation}
and it is easy to see that this gives a diffeomorphism only near the origin.
\\[0.2cm]
\textbf{(II)} \, Let us calculate the $C^\infty$ normal form of Eq.(\ref{4-2}).
The first term $\mathcal{P}_K(R_1)$ of the normal form is given by using Eq.(\ref{hyp3-1}) as
\begin{equation}
\mathcal{P}_K(R_1) (y_1, y_2) = \lim_{t\to \infty} \frac{1}{t} \int^t \! \left(
\begin{array}{@{\,}cc@{\,}}
e^{-is} & 0 \\
0 & e^{is}
\end{array}
\right) \left(
\begin{array}{@{\,}c@{\,}}
\sin (e^{is} y_1 + e^{-is} y_2) \\
\sin (e^{is} y_1 + e^{-is} y_2) 
\end{array} 
\right) ds .
\label{4-7}
\end{equation}
Thus the first order normal form is given by
\begin{equation}
\frac{d}{dt}\left(
\begin{array}{@{\,}c@{\,}}
y_1 \\
y_2
\end{array}
\right) = \left(
\begin{array}{@{\,}c@{\,}}
iy_1 \\
-iy_2
\end{array}
\right) + \frac{\varepsilon }{2\pi} \left(
\begin{array}{@{\,}c@{\,}}
\int^{2\pi}_{0} \! e^{-it} \sin (e^{it}y_1 + e^{-it}y_2) dt\\
 \int^{2\pi}_{0} \! e^{it} \sin (e^{it}y_1 + e^{-it}y_2) dt
\end{array}
\right) .
\label{4-8}
\end{equation}
Putting $y_1 = re^{i\theta },\, y_2 = re^{-i\theta }$ yields
\begin{equation}
\left\{ \begin{array}{l}
\displaystyle \dot{r} = \frac{\varepsilon }{2\pi} \int^{2\pi}_{0} \! \cos t \cdot \sin (2r \cos t) dt = \varepsilon J_1 (2r),\\
\displaystyle \dot{\theta } = 1 + \frac{\varepsilon }{2\pi r}\int^{2\pi}_{0} \! \sin t \cdot \sin (2r \cos t) dt = 1,
\end{array} \right.
\label{4-9}
\end{equation}
where $J_n(r)$ is the Bessel function of the first kind defined as the solution of the equation
$r^2 x'' + rx' + (r^2 - n^2)x = 0$.
By Eq.(\ref{hyp3-2}), it is easy to verify that the first order near identity transformation
is periodic in $y_1$ and $y_2$ although we can not calculate the indefinite integral in Eq.(\ref{hyp3-2}) explicitly.
Thus there exists a positive number $\varepsilon _0$ such that if $0 < \varepsilon  < \varepsilon _0$,
the near identity transformation is a diffeomorphism on $\mathbf{R}^2$.
Since $J_1 (2r)$ has infinitely many zeros, Thm.3.9 proves that the original system (\ref{4-1})
has infinitely many periodic orbits.
\\[0.2cm]
\textbf{Example \thedef.} Consider the system on $\mathbf{R}^2$ of the form
\begin{equation}
\left\{ \begin{array}{l}
\dot{x}_1 = x_2 + 2\varepsilon g(x_1),  \\
\dot{x}_2 = -x_1 , \\
\end{array} \right.
\label{4-10}
\end{equation}
where the function $g(x)$ is defined by
\begin{equation}
g(x) = \left\{ \begin{array}{l}
x, \quad x\in [ 2n, 2n+1 ),  \\
-x, \quad x\in [ 2n+1, 2n+2 ),  \\
\end{array} \right.
\label{4-11}
\end{equation}
for $n=0,1,2,\cdots $ and $g(x) = - g(-x)$ (see Fig.1 (a)).

\begin{figure}[h]
\begin{center}
\includegraphics[scale=1.0]{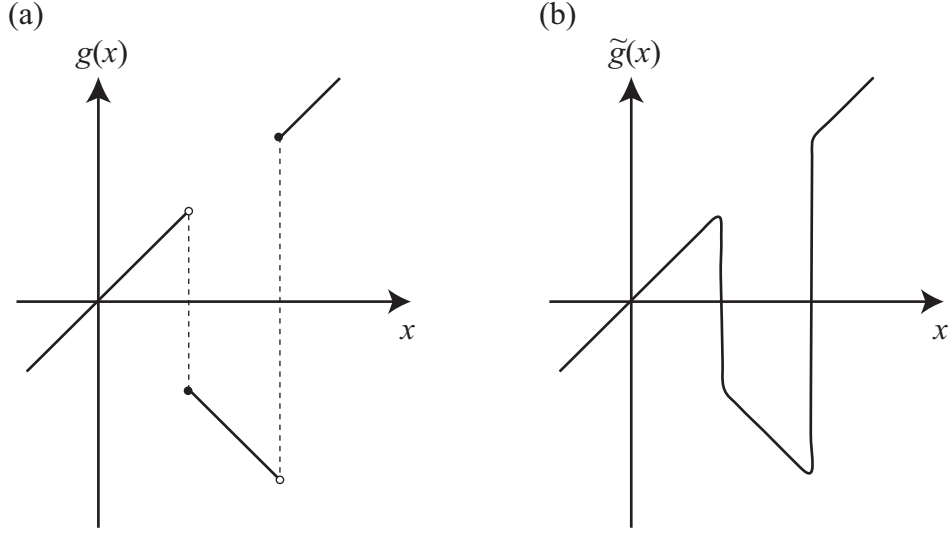}
\end{center}
\caption{The graphs of the functions $g(x)$ and $\widetilde{g}(x)$.}
\end{figure}

We add to Eq.(\ref{4-10}) a small perturbation whose support is included in 
sufficiently small intervals $(n- \delta , n+ \delta ),\, n \in \mathbf{Z}$ so that
the resultant system
\begin{equation}
\left\{ \begin{array}{l}
\dot{x}_1 = x_2 + 2\varepsilon \widetilde{g}(x_1),  \\
\dot{x}_2 = -x_1 , \\
\end{array} \right.
\label{4-12}
\end{equation}
is of $C^\infty$ class (see Fig.1 (b)).
Like as Example 4.1, the first order $C^\infty$ normal form of this system written in the polar coordinates
is given by
\begin{equation}
\left\{ \begin{array}{l}
\displaystyle \dot{r} = \frac{\varepsilon }{2\pi}\int^{2\pi}_{0} \! \cos t \cdot \widetilde{g}(2r \cos t) dt
   := \frac{\varepsilon }{2\pi} R(r),   \\
\displaystyle \dot{\theta } = 1 + \frac{\varepsilon }{2\pi}\int^{2\pi}_{0} \! \sin t \cdot \widetilde{g}(2r \cos t) = 1. \\
\end{array} \right.
\label{4-13}
\end{equation}
On the outside of the support of the perturbation, the function $R(r)$ is given by
\begin{equation}
R(r) = 
\left\{ \begin{array}{ll}
2\pi r, & r\in ( 2n + \delta , 2n +1 - \delta), \\
-2\pi r, & r\in ( 2n+1 + \delta , 2n +2 - \delta). \\
\end{array} \right.
\end{equation}
By the intermediate value theorem, $R(r)$ has zeros near $r = n\in \mathbf{Z}$.
In particular, fixed points near $r = 2n+1$ is attracting.
This and Thm.3.9 prove that Eq.(\ref{4-12}) has stable periodic orbits.
Therefore Eq.(\ref{4-10}) has attracting invariant sets near the stable periodic orbits of Eq.(\ref{4-12}).

If we apply the polynomial normal form to Eq.(\ref{4-12}) after expanding $\widetilde{g}(x)$ at the origin,
we obtain the normal form $\dot{r} = \varepsilon r$, which is valid on a small neighborhood of the origin.

\appendix

\section{Appendix}

In this appendix, we derive Eq.(\ref{3-33}) from Eq.(\ref{3-32}).
By integrating by parts, Eq.(\ref{3-32}) is calculated as
\begin{eqnarray*}
x_2 &=& e^{At}h^{(2)}(y) + e^{At}\int^t_{0} \!e^{-As} \left( 
         \frac{\partial g_1}{\partial x}(e^{As}y) \mathcal{Q}(g_{1I})(e^{As}y) + g_2 (e^{As}y) \right) ds \nonumber \\
 & & + e^{At}\!\!\int^t_{0} \!  e^{-As} \frac{\partial g_1}{\partial x}(e^{As}y) g_{1K}(e^{As}y) ds \cdot t 
     - e^{At}\!\!\int^t_{0} \!\!ds \! \int^s_{0} \!\!
         e^{-As'} \frac{\partial g_1}{\partial x}(e^{As'}y)g_{1K}(e^{As'}y) ds' \nonumber \\
&=& e^{At}h^{(2)}(y) + e^{At}\int^t_{0} \!e^{-As} \left( 
         \frac{\partial g_1}{\partial x}(e^{As}y) \mathcal{Q}(g_{1I})(e^{As}y) + g_2 (e^{As}y) \right) ds \nonumber \\
 & & +  e^{At}\!\!\int^t_{0} \!  e^{-As} \frac{\partial g_{1K}}{\partial x}(e^{As}y) g_{1K}(e^{As}y) ds \cdot t 
     + e^{At}\!\!\int^t_{0} \!  e^{-As} \frac{\partial g_{1I}}{\partial x}(e^{As}y) g_{1K}(e^{As}y) ds \cdot t \nonumber \\
& & - e^{At}\!\!\int^t_{0} \!\!ds \! \int^s_{0} \!\!
         e^{-As'} \frac{\partial g_{1K}}{\partial x}(e^{As'}y)g_{1K}(e^{As'}y) ds'
   - e^{At}\!\!\int^t_{0} \!\!ds \! \int^s_{0} \!\!
         e^{-As'} \frac{\partial g_{1I}}{\partial x}(e^{As'}y)g_{1K}(e^{As'}y) ds'.
\end{eqnarray*}
Since $Dg_{1K} \cdot g_{1K}\in V_K$ and $Dg_{1I} \cdot g_{1K}\in V_I$ by Props.3.3 and 3.4, we obtain
\begin{eqnarray*}
x_2 &=& e^{At}h^{(2)}(y) + e^{At}\int^t_{0} \!e^{-As} \left( 
         \frac{\partial g_1}{\partial x}(e^{As}y) \mathcal{Q}(g_{1I})(e^{As}y) + g_2 (e^{As}y) \right) ds \nonumber \\
 & & +  e^{At} \frac{\partial g_{1K}}{\partial x}(y) g_{1K}(y) t^2  
     + \mathcal{Q} \left( \frac{\partial g_{1I}}{\partial x} g_{1K}\right) (e^{At}y) t
     - e^{At}\mathcal{Q} \left( \frac{\partial g_{1I}}{\partial x} g_{1K}\right) (y) t  \nonumber \\
& & - e^{At} \int^t_{0}\! \frac{\partial g_{1K}}{\partial x}(y)g_{1K}(y) s ds
    - e^{At}\!\int^t_{0} \!  
      \left( e^{-As} \mathcal{Q} \left( \frac{\partial g_{1I}}{\partial x} g_{1K}\right) (e^{As}y)
           - \mathcal{Q} \left( \frac{\partial g_{1I}}{\partial x} g_{1K}\right) (y) \right) ds \nonumber \\
&=& e^{At}h^{(2)}(y) + e^{At}\int^t_{0} \!e^{-As} \left( 
         \frac{\partial g_1}{\partial x} \mathcal{Q}(g_{1I})+ g_2 
         - \frac{\partial \mathcal{Q}(g_{1I})}{\partial x} g_{1K} \right) (e^{As}y) ds \nonumber \\
 & & + \frac{1}{2} e^{At} \frac{\partial g_{1K}}{\partial x}(y) g_{1K}(y) t^2  
      + \frac{\partial \mathcal{Q}(g_{1I})}{\partial x}(e^{At}y) g_{1K} (e^{At}y) t.
\end{eqnarray*}
Since $R_2$ is defined by Eq.(\ref{3-34}), the above is rewritten as
\begin{eqnarray*}
x_2 &=& e^{At}h^{(2)}(y) + e^{At}\int^t_{0} \!e^{-As} \mathcal{P}_I(R_2) (e^{As}y) ds
       + e^{At}\int^t_{0} \!e^{-As} \mathcal{P}_K(R_2) (e^{As}y) ds \nonumber \\
 & & + \frac{1}{2} e^{At} \frac{\partial g_{1K}}{\partial x}(y) g_{1K}(y) t^2  
      + \frac{\partial \mathcal{Q}(g_{1I})}{\partial x}(e^{At}y) g_{1K} (e^{At}y) t \nonumber \\
&=& e^{At}h^{(2)}(y) + \mathcal{Q}\mathcal{P}_I (R_2)(e^{At}y) - e^{At} \mathcal{Q}\mathcal{P}_I(R_2)(y)
       + e^{At} \mathcal{P}_K(R_2) (y) t \nonumber \\
 & & + \frac{1}{2} e^{At} \frac{\partial g_{1K}}{\partial x}(y) g_{1K}(y) t^2  
      + \frac{\partial \mathcal{Q}(g_{1I})}{\partial x}(e^{At}y) g_{1K} (e^{At}y) t.
\end{eqnarray*}
Putting $h^{(2)} = \mathcal{Q}\mathcal{P}_I (R_2)$, we obtain Eq.(\ref{3-33}).


\end{document}